\documentclass[11pt,oneside,french,english]{amsart}
\usepackage[T1]{fontenc}
\usepackage[latin9]{inputenc}
\usepackage[letterpaper]{geometry}
\usepackage{amsthm}
\usepackage{amssymb}

\makeatletter

\newcommand{\noun}[1]{\textsc{#1}}

\numberwithin{equation}{section}
\numberwithin{figure}{section}

\makeatother

\usepackage{babel}
\addto\extrasfrench{\providecommand{\fg}{\ifdim\lastskip>\z@\unskip\fi~\frqq}}

\begin{document}
\selectlanguage{french}%
\textbf{\hfill{}} \textit{Paru à }\foreignlanguage{english}{\textit{Comptes
Rendus de l'Académie des Sciences de Paris, Ser. I 347 (2009) 785\textendash{}790.}}

\selectlanguage{english}%
\vspace{0.5cm}

\begin{center}
\noun{\Large Une amélioration d'un résultat de E}{\Large . }\noun{\Large B}{\Large .}\noun{\Large{}
Davies et B. Simon}
\par\end{center}{\Large \par}

\noun{\large \vspace{0.2cm}
}{\large \par}

\begin{center}
\textit{\noun{\small Rachid Zarouf}}
\par\end{center}{\small \par}

\noun{\large \vspace{0.2cm}
}{\large \par}

\begin{flushleft}
\textbf{\small Résumé}
\par\end{flushleft}{\small \par}

\begin{flushleft}
{\small E. B. Davies et B. Simon ont montré (entre autres résultats)
la chose suivante : soit $T$ , une matrice $n\times n$ telle que
son spectre $\sigma(T)$ soit inclus dans le disque $\mathbb{D}=\left\{ z\in\mathbb{C}:\,\vert z\vert<1\right\} $
et soit $C=sup_{n\geq0}\left\Vert T^{n}\right\Vert _{E\rightarrow E}$
, ($E$ étant $\mathbb{C}^{n}$ muni de la norme euclidienne $\left|.\right|$
). Alors $\left\Vert R(1,\, T)\right\Vert _{E\rightarrow E}\leq C\left(3n/dist(1,\,\sigma(T))\right)^{3/2}$
où $R(\lambda,\, T)$ désigne la résolvante de $T$ prise au point
$\lambda.$ Nous améliorons ici cette dernière inégalité à travers
le résultat suivant : sous les mêmes conditions (portant sur la matrice
$T$) et si $E$ est cette fois $\mathbb{C}^{n}$ muni d'une norme
quelconque $\left|.\right|$, alors pour tout $\lambda\notin\sigma(A)$
tel que $\vert\lambda\vert\geq1$, on a $\left\Vert R(\lambda,\, T)\right\Vert \leq C\left(5\pi/3+2\sqrt{2}\right)n^{3/2}/dist\left(\lambda,\,\sigma\right).$}
\par\end{flushleft}{\small \par}

\noun{\small \vspace{0.1cm}
}{\small \par}

\selectlanguage{french}%
\begin{center}
\noun{\small Sharpening a result by E.B. Davies and B. Simon}\foreignlanguage{english}{\noun{\small \vspace{0.1cm}
}}
\par\end{center}{\small \par}

\selectlanguage{english}%
\begin{flushleft}
\textbf{\small Abstract}
\par\end{flushleft}{\small \par}

\begin{flushleft}
{\small E. B. Davies et B. Simon have shown (among other things) the
following result}\noun{\small :}{\small{} if $T$ is an $n\times n$
matrix such that its spectrum $\sigma(T)$ is included in the open
unit disc $\mathbb{D}=\left\{ z\in\mathbb{C}:\,\vert z\vert<1\right\} $
and if $C=sup_{k\geq0}\left\Vert T^{k}\right\Vert _{E\rightarrow E}$,
where $E$ stands for $\mathbb{C}^{n}$ endowed with the euclidean
norm $\left|.\right|_{2}$ , then $\left\Vert R(1,\, T)\right\Vert _{E\rightarrow E}\leq C\left(3n/dist(1,\,\sigma(T))\right)^{3/2}$
where $R(\lambda,\, T)$ stands for the resolvent of $T$ at point
$\lambda$. Here, we improve this inequality showing that under the
same hypotheses (on the matrix $T$) and if this time $E$ is $\mathbb{C}^{n}$
endowed with certain norm $\left|.\right|$ , then$\left\Vert R(\lambda,\, T)\right\Vert \leq C\left(5\pi/3+2\sqrt{2}\right)n^{3/2}/dist\left(\lambda,\,\sigma\right),$
for all $\lambda\notin\sigma(T)$ such that $\vert\lambda\vert\geq1$. }
\par\end{flushleft}{\small \par}

\noun{\small \vspace{0.2cm}
}{\small \par}

\begin{flushleft}
Pour $C\geq1$ et $n\in\mathbb{N}^{\star}$ on pose \[
K_{n}(C)=sup\left\Vert R(\lambda,\, T)\right\Vert dist\left(\lambda,\,\sigma(T)\right),\]

\par\end{flushleft}

\noun{\large \vspace{0.05cm}
}{\large \par}

\begin{flushleft}
où la borne supérieure est prise sur l'ensemble des $\lambda\in\mathbb{C}$
tels que $\left|\lambda\right|\geq1$ et sur l'ensemble des opérateurs
$T:\, E\rightarrow E$ avec $E=\left(\mathbb{C}^{n},\,\left|.\right|\right)$
et vérifiant $\forall k\in\mathbb{N},$ $\left\Vert T^{k}\right\Vert _{E\rightarrow E}\leq C$.
\par\end{flushleft}

Le but de cette note est de démontrer le théorème suivant.

\begin{flushleft}
\textbf{Théorème. }\textit{(i) Pour tout $n\in\mathbb{N}^{\star}$
et pour tout $C\geq1,$ on a}
\par\end{flushleft}

\textit{\[
K_{n}(C)\leq C\left(5\pi/3+2\sqrt{2}\right)n^{3/2}.\]
}

\textit{(ii) De plus, pour tout $C\geq1$ on a \[
limsup_{n\rightarrow\infty}n^{-\frac{3}{2}}K_{n}(C)\leq5C\pi/3.\]
}

\begin{flushleft}
\textbf{Commentaires.}
\par\end{flushleft}

(1) Ce théorème est un résultat de nature {}``numérique'' dans la
mesure où il s'agit d'une estimation du type $\left\Vert R(\lambda,\, T)\right\Vert \leq Kdist\left(\lambda,\,\sigma(T)\right),$
où la question est d'évaluer la taille de la constante $K=K_{n}(C)$
en fonction des paramètres dont elle dépend, à savoir $n$ et $C.$

(2) Le résultat principal de E. B. Davies et B. Simon est le suivant,
voir {[}2{]}: Soit $K_{n}$ la même borne supérieure que celle donnant
$K_{n}(C)$ en restreignant la condition $(C_{1}):\:${[}\textit{$T:\, E\rightarrow E$
avec $E=\left(\mathbb{C}^{n},\,\left|.\right|\right)$ et vérifiant
$\forall k\in\mathbb{N},$ $\left\Vert T^{k}\right\Vert _{E\rightarrow E}\leq C$}{]}
par $(C_{2}):$ {[}\textit{T est une contraction d'un espace de Hilbert}{]}.
Alors le facteur $n^{\frac{3}{2}}$ {}``devient'' $n$ et $K_{n}=cotan(\pi/4n).$ 

(3) Dans ce dernier cas (où $T$ est une contraction d'un espace de
Hilbert), la méthode appliquée ci-dessous pour montrer le Théorème
faisant l'objet de cette note, donne $K_{n}\leq an$ où $a=(1+max_{\lambda\in\sigma(T)}\left|\lambda\right|).$
En particulier, pour $r=max_{\lambda\in\sigma(T)}\left|\lambda\right|<4/\pi-1$,
cette majoration est plus précise que {[}2{]}.

(4) En ce qui concerne l'exactitude du véritable ordre de croissance
de la constante $K_{n}(C)$, on sait pour l'instant que $K_{n}(C)/n\geq K_{n}/n\geq b$
où $b=(2+\sqrt{3})/3$, voir {[}2{]} p.4.

(5) L'hypothèse du théorème entraîne trivialement que $\left\Vert R(\lambda,T)\right\Vert \leq C\left(\left|\lambda\right|-1\right)^{-1}$,
$\left|\lambda\right|>1.$ Ce théorème peut donc être vu comme un
analogue unilatéral du {}``Lemme de Domar'' bien connu (voir {[}1{]},
{[}6{]}): si $\sigma\subset\mathbb{D}$ et $u$ une fonction sous-harmonique
dans $\mathbb{C}\setminus\sigma$ telle que pour tout $\lambda\in\mathbb{C}\setminus\sigma,$
$u(\lambda)\leq Cmax\left\{ \left|\left|\lambda\right|-1\right|^{-1},\, dist\left(\lambda,\,\sigma\right)^{-1}\right\} $,
alors pour tout $\lambda\in\mathbb{C}\setminus\sigma,$ tel que $\left|\lambda\right|\geq1/2$,
$u(\lambda)\leq447Cdist\left(\lambda,\,\sigma\right)^{-1}.$ Une version
aussi générale pour des estimations unilatérales $\left(\left|\lambda\right|\geq1\right)$,
n'est pas vraie (exemple: $u(\lambda)=\left\Vert R(\lambda,M_{\theta})\right\Vert $,
où $M_{\theta}$ est l'opérateur modèle sur $K_{\theta}=H^{2}\Theta\theta H^{2}$,
$\theta=exp\left(\frac{z+1}{z-1}\right),$ voir {[}4{]}), mais notre
résultat montre qu'elle est correcte pour des résolvantes de matrices
de taille $n$ (avec une constante dépendant de $n$).

\noun{\small \vspace{0.1cm}
}{\small \par}

Nous avons recours en premier lieu au lemme suivant, de type principe
du maximum.

\begin{flushleft}
\textbf{Lemme.}\textit{ Soient $C\geq1$, $A>0$ tels que pour tout
opérateur $T$ agissant sur $\left(\mathbb{C}^{n},\,\left|.\right|\right)$
et de spectre $\sigma(T)$, la condition suivante soit réalisée: \[
\left\{ \begin{array}{c}
sup_{k\geq0}\left\Vert T^{k}\right\Vert \leq C\\
\sigma(T)\subset\mathbb{D}\end{array}\right.\Longrightarrow\left[\forall\lambda_{\star}\: tel\, que\:\left|\lambda_{\star}\right|=1,\: dist\left(\lambda_{\star},\,\sigma(T)\right)\left\Vert R\left(\lambda_{\star},T\right)\right\Vert \leq A\right],\]
alors, \[
K_{n}(C)\leq A.\]
}
\par\end{flushleft}

\begin{flushleft}
\textbf{Preuve.} Soit $\lambda$ tel que $\left|\lambda\right|>1$.
$\lambda$ peut alors s'écrire $\lambda=\rho\lambda_{\star}$ avec
$\rho>1$ et $\left|\lambda_{\star}\right|=1.$ On pose $T_{\star}=\frac{1}{\rho}T$.
Dans ces conditions, $sup_{k\geq0}\left\Vert T_{\star}^{k}\right\Vert \leq C$
et $\sigma(T_{\star})=\frac{1}{\rho}\sigma(T)\subset\mathbb{D}$.
Par conséquent, on a $dist\left(\lambda_{\star},\,\sigma\left(T_{\star}\right)\right)\left\Vert R\left(\lambda_{\star},T_{\star}\right)\right\Vert \leq A,$
ce que l'ont peut encore écrire $\rho dist\left(\lambda_{\star},\,\sigma\left(T_{\star}\right)\right)\left\Vert \rho^{-1}R\left(\lambda_{\star},T_{\star}\right)\right\Vert \leq A.$
Il suffit maintenant de remarquer que $\rho dist\left(\lambda_{\star},\,\sigma\left(T_{\star}\right)\right)=dist\left(\lambda,\,\sigma(T)\right)$
et $\rho^{-1}R\left(\lambda_{\star},T_{\star}\right)=R(\lambda,T).$ 
\par\end{flushleft}

\begin{flushright}
$\square$
\par\end{flushright}

\begin{flushleft}
\textbf{Preuve du Théorème.} Soient $T$ une matrice de taille $n$
vérifiant la condition $(C_{1})$ et $\sigma=\sigma(T)=\left\{ \lambda_{1},\,\lambda_{2},\,...,\,\lambda_{n}\right\} $
son spectre (les $\lambda_{j}$ étant comptés avec leur multiplicité).
On défini le produit de Blaschke $B=\Pi_{k=1}^{n}b_{\lambda_{i}},$
où pour tout $i=1..n$, $b_{\lambda_{i}}=\frac{\lambda_{i}-z}{1-\overline{\lambda_{i}}z}$
. Tout d'abord,
\par\end{flushleft}

\textit{\[
\left\Vert R(\lambda,T)\right\Vert \leq C\left\Vert \frac{1}{\lambda-z}\right\Vert _{W/BW},\]
}

\noun{\large \vspace{0.05cm}
}{\large \par}

\begin{flushleft}
(voir {[}3{]} Théorème 3.24, p.31), où $W$ est l'algèbre de Wiener
des séries de Taylor absolument convergentes, $W=\left\{ f=\sum_{k\geq0}\hat{f}(k)z^{k}:\:\left\Vert f\right\Vert _{W}=\sum_{k\geq0}\left|\hat{f}(k)\right|<\infty\right\} $
et\textit{ }
\par\end{flushleft}

\textit{\[
\left\Vert \frac{1}{\lambda-z}\right\Vert _{W/BW}=inf\left\{ \left\Vert f\right\Vert _{W}:\, f\left(\lambda_{j}\right)=\frac{1}{\lambda-\lambda_{j}},\, j=1..n\right\} .\]
}

\noun{\large \vspace{0.05cm}
}{\large \par}

\begin{flushleft}
On suppose dans un premier temps que $\vert\lambda\vert>1.$ Soit
$P_{B}$ la projection orthogonale de l'espace de Hardy $H^{2}$ sur
$K_{B}=H^{2}\Theta BH^{2}$. La fonction $f=P_{B}(\frac{1}{\lambda}k_{1/\bar{\lambda}})$
verifie bien $f-\frac{1}{\lambda-z}\in BW,\,\forall\, j=1..n$. En
particulier, on a \[
\left\Vert \frac{1}{\lambda-z}\right\Vert _{W/BW}\leq\left\Vert \frac{1}{\lambda}P_{B}k_{1/\bar{\lambda}}\right\Vert _{W}.\]
Mais on sait que
\par\end{flushleft}

\[
P_{B}k_{1/\bar{\lambda}}=\sum_{k=1}^{n}\left(k_{1/\bar{\lambda}},\, e_{k}\right)_{H^{2}}e_{k}\]
où la famille $\left(e_{k}\right)_{k=1}^{n}$ (appelée base de Malmquist
relative à $\sigma$, voir {[}5{]} p.117) définie par, \[
e_{1}=\left(1-\vert\lambda_{1}\vert^{2}\right)^{\frac{1}{2}}f_{1},\: e_{k}=\left(1-\left|\lambda_{k}\right|^{2}\right)^{\frac{1}{2}}\left(\Pi_{j=1}^{k-1}b_{\lambda_{j}}\right)f_{k}=(f_{k}/\left\Vert f_{k}\right\Vert _{2})\Pi_{j=1}^{k-1}b_{\lambda_{j}},\, k\geq2,\]
où $f_{k}(z)=\frac{1}{1-\overline{\lambda_{k}}z}.$ Du coup,

\[
P_{B}k_{1/\bar{\lambda}}=\sum_{k=1}^{n}\overline{e_{k}\left(1/\bar{\lambda}\right)}e_{k}.\]
Nous allons maintenant appliquer l'inégalité de Hardy $\left\Vert f\right\Vert _{W}\leq\pi\left\Vert f'\right\Vert _{H^{1}}+\left|f(0)\right|$,
(voir N. Nikolski, {[}4{]} p. 370 8.7.4 -(c)) à $P_{B}k_{1/\bar{\lambda}}$
en profitant du fait remarquable que pour $k=2..n$\[
e_{k}^{'}=\sum_{i=1}^{k-1}\frac{b_{\lambda_{i}}^{'}}{b_{\lambda_{i}}}e_{k}+\overline{\lambda_{k}}\frac{1}{\left(1-\overline{\lambda_{k}}z\right)}e_{k}.\]
On trouve alors

\[
\left(P_{B}k_{1/\bar{\lambda}}\right)^{'}=\left(k_{1/\bar{\lambda}},\, e_{1}\right)_{H^{2}}\frac{\bar{\lambda}_{1}}{\left(1-\overline{\lambda_{1}}z\right)}e_{1}+\sum_{i=1}^{n}\frac{b_{\lambda_{i}}^{'}}{b_{\lambda_{i}}}\sum_{k=i+1}^{n-1}\left(k_{1/\bar{\lambda}},\, e_{k}\right)_{H^{2}}e_{k}+\sum_{k=2}^{n}\left(k_{1/\bar{\lambda}},\, e_{k}\right)_{H^{2}}\overline{\lambda_{k}}\frac{1}{\left(1-\overline{\lambda_{k}}z\right)}e_{k}.\]
Comme $k_{1/\bar{\lambda}}$ est le noyau reproduisant de $H^{2}$
associé au point $1/\bar{\lambda}\in\mathbb{D}$, on trouve $\left(e_{k},\, k_{1/\bar{\lambda}}\right)_{H^{2}}=e_{k}(1/\bar{\lambda}),$
et donc\[
\left(P_{B}k_{1/\bar{\lambda}}\right)^{'}=\overline{e_{1}\left(1/\bar{\lambda}\right)}\frac{\bar{\lambda}_{1}}{\left(1-\overline{\lambda_{1}}z\right)}e_{1}+\sum_{i=1}^{n}\frac{b_{\lambda_{i}}^{'}}{b_{\lambda_{i}}}\sum_{k=i+1}^{n-1}\overline{e_{k}\left(1/\bar{\lambda}\right)}e_{k}+\sum_{k=2}^{n}\overline{e_{k}\left(1/\bar{\lambda}\right)}\overline{\lambda_{k}}\frac{1}{\left(1-\overline{\lambda_{k}}z\right)}e_{k}.\]
Maintenant,

\[
\left\Vert e_{1}\left(1/\bar{\lambda}\right)\frac{\lambda_{1}}{\left(1-\overline{\lambda_{1}}z\right)}e_{1}\right\Vert _{H^{1}}\leq\left|e_{1}\left(1/\bar{\lambda}\right)\right|\left\Vert \frac{\lambda_{1}}{\left(1-\overline{\lambda_{1}}z\right)}\right\Vert _{H^{2}}\left\Vert e_{1}\right\Vert _{H^{2}}\leq\left|\lambda\right|\frac{1}{dist\left(\lambda,\,\sigma\right)}\]
en utilisant à la fois l'inégalité de Cauchy-Schwarz et le fait que
$e_{1}$ est de norme 1 dans $H^{2}$. Par la même raison (la famille
$\left(e_{k}\right)_{k=1}^{n}$ est orthonormale dans $H^{2}$), on
trouve\[
\left\Vert \sum_{k=2}^{n}\overline{\lambda_{k}}\overline{e_{k}\left(1/\bar{\lambda}\right)}\frac{1}{\left(1-\overline{\lambda_{k}}z\right)}e_{k}\right\Vert _{H^{1}}\leq\sum_{k=2}^{n}\left|e_{k}\left(1/\bar{\lambda}\right)\right|\left\Vert \lambda_{k}\frac{1}{\left(1-\overline{\lambda_{k}}z\right)}\right\Vert _{H^{2}}\left\Vert e_{k}\right\Vert _{H^{2}}\leq\]
\[
\leq\sum_{k=2}^{n}\left|\frac{\left(1-\left|\lambda_{k}\right|^{2}\right)^{\frac{1}{2}}}{1-\overline{\lambda_{k}}/\overline{\lambda}}\right|\frac{1}{\sqrt{1-\left|\lambda_{k}\right|^{2}}}\leq\left|\lambda\right|\frac{(n-1)}{dist\left(\lambda,\,\sigma\right)}.\]
Finalement,\[
\left\Vert \sum_{i=1}^{n-1}\frac{b_{\lambda_{i}}^{'}}{b_{\lambda_{i}}}\sum_{k=i+1}^{n}\overline{e_{k}\left(1/\bar{\lambda}\right)}e_{k}\right\Vert _{H^{1}}\leq\sum_{i=1}^{n-1}\left\Vert \frac{b_{\lambda_{i}}^{'}}{b_{\lambda_{i}}}\right\Vert _{L^{2}}\left(\sum_{k=i+1}^{n}\left|e_{k}\left(1/\bar{\lambda}\right)\right|^{2}\right)^{\frac{1}{2}}.\]
Comme en outre on a $b_{\lambda_{i}}^{'}/b_{\lambda_{i}}=1/\left(\lambda_{i}-z\right)+\overline{\lambda_{i}}/\left(1-\overline{\lambda_{i}}z\right),$
on en déduit que $\left\Vert b_{\lambda_{i}}^{'}/b_{\lambda_{i}}\right\Vert _{L^{2}}\leq2/\sqrt{1-\left|\lambda_{i}\right|^{2}}$,
et que\[
\left\Vert \sum_{i=1}^{n-1}\frac{b_{\lambda_{i}}^{'}}{b_{\lambda_{i}}}\sum_{k=i+1}^{n}\overline{e_{k}\left(1/\bar{\lambda}\right)}e_{k}\right\Vert _{H^{1}}\leq2\sum_{i=1}^{n-1}\frac{1}{\left(1-\left|\lambda_{i}\right|^{2}\right)^{\frac{1}{2}}}\left(\sum_{k=i+1}^{n}\left|\frac{\left(1-\left|\lambda_{k}\right|^{2}\right)^{}}{\left(1-\overline{\lambda_{k}}/\bar{\lambda}\right)^{2}}\right|\right)^{\frac{1}{2}}.\]
Maintenant, sans perte de généralité on peut supposer que la suite
$\left(\left|\lambda_{i}\right|\right)_{i=1}^{n}$ est croissante
(quitte à réordonner la séquence $\sigma$). Dans ce cas, pour $k\geq i+1>i$
on a $1-\left|\lambda_{k}\right|^{2}\leq1-\left|\lambda_{i}\right|^{2}$
ce qui donne\[
\left\Vert \sum_{i=1}^{n-1}\frac{b_{\lambda_{i}}^{'}}{b_{\lambda_{i}}}\sum_{k=i+1}^{n}\overline{e_{k}\left(1/\bar{\lambda}\right)}e_{k}\right\Vert _{H^{1}}\leq2\sum_{i=1}^{n-1}\left(\sum_{k=i+1}^{n}\left|\frac{1}{\left(1-\overline{\lambda_{k}}/\bar{\lambda}\right)^{2}}\right|\right)^{\frac{1}{2}}\leq\]
\[
\leq2\frac{\left|\lambda_{}\right|}{dist\left(\lambda,\,\sigma\right)}\sum_{i=1}^{n-1}\left(\sum_{k=i+1}^{n}1\right)^{\frac{1}{2}}\leq\frac{4}{3}\left|\lambda_{}\right|\frac{1}{dist\left(\lambda,\,\sigma\right)}\left(n^{\frac{3}{2}}-1\right),\]
puisque\textbf{\textit{ }}$\sum_{i=1}^{n-1}\sqrt{j}\leq\int_{1}^{n}\sqrt{x}dx$.
Finalement,

\[
\left\Vert \left(\frac{1}{\lambda}P_{B}k_{1/\bar{\lambda}}\right)^{'}\right\Vert _{H^{1}}\leq\frac{1}{dist\left(\lambda,\,\sigma\right)}+\frac{(n-1)}{dist\left(\lambda,\,\sigma\right)}+\frac{4}{3}\frac{n^{\frac{3}{2}}-1}{dist\left(\lambda,\,\sigma\right)}=\]
\[
=\frac{1}{dist\left(\lambda,\,\sigma\right)}\left(-\frac{4}{3}+n+\frac{4}{3}n^{\frac{3}{2}}\right)\leq\frac{5}{3}\frac{n^{\frac{3}{2}}}{dist\left(\lambda,\,\sigma\right)},\]
la dernière inégalité reposant sur le fait que pour tout $x\geq0,$
$\frac{1}{3}x^{\frac{3}{2}}-x+\frac{4}{3}\geq0$. Ceci donne\[
\left\Vert \frac{1}{\lambda}P_{B}k_{1/\bar{\lambda}}\right\Vert _{W}\leq\frac{5}{3}\pi\frac{n^{\frac{3}{2}}}{dist\left(\lambda,\,\sigma\right)}+\left|\frac{1}{\lambda}\right|\sum_{k=1}^{n}\left|e_{k}\left(1/\bar{\lambda}\right)\right|\left|e_{k}(0)\right|\leq\frac{5}{3}\pi\frac{n^{\frac{3}{2}}}{dist\left(\lambda,\,\sigma\right)}+2n.\]
En particulier, (ii) est démontré. Pour résumer, on a pour $n\geq1$,\[
\left\Vert R(\lambda,T)\right\Vert \leq C\left(5\pi/3+\frac{2}{\sqrt{n}}dist\left(\lambda,\,\sigma\right)\right)\frac{n^{\frac{3}{2}}}{dist\left(\lambda,\,\sigma\right)}.\]
Faisons maintenant tendre radialement $\lambda$ vers sa projection
$\lambda_{\star}$ sur le tore $\mathbb{T}$ et remarquons qu'alors,
(puisque $dist\left(\lambda_{\star},\,\sigma\right)\leq2$), \[
\left\Vert R\left(\lambda_{\star},T\right)\right\Vert \leq C\left(5\pi/3+2\sqrt{2}\right)\frac{n^{\frac{3}{2}}}{dist\left(\lambda_{\star},\,\sigma\right)}.\]
Il reste alors à appliquer le lemme avec $A=\left(5\pi/3+2\sqrt{2}\right)n^{\frac{3}{2}}$
pour achever la preuve de (i).

\begin{flushright}
$\square$
\par\end{flushright}

\begin{flushleft}
\textbf{Remerciements}
\par\end{flushleft}

\vspace{0.1cm}

Je tiens à remercier infiniment le Professeur Nikolai Nikolski pour
ses conseils ô combien précieux.

\noun{\small \vspace{0.1cm}
}{\small \par}

\[
\begin{array}{c}
\mbox{Equipe\:\ d'Analyse\:\ et\:\ Géométrie,}\\
\mbox{Institut\:\ de\:\ Mathématiques\:\ de\:\ Bordeaux,}\\
\mbox{\mbox{Université\:\ Bordeaux,\:351\:\ Cours\:\ de\:\ la\:\ Libération,\,33405\:\ Talence,\:\ France}.}\\
\mbox{E-mail\:\ address:\:\ rzarouf@math.u-bordeaux1.fr}\end{array}\]


\begin{thebibliography}{6}
\bibitem[1]{key-1} Y. Domar, \textit{On the existence of a largest
subharmonic minorant of a given function}, On the existence of a largest
subharmonic minorant of a given function, 3 (1958), 429\textendash{}440

\bibitem[2]{key-1} E. B. Davies and B. Simon, \textit{Eigenvalue
estimates for non-normal matrices and the zeros of random orthogonal
polynomials on the unit circle}, J. Approx. Theory 141-2, (2006),
189\textendash{}213.

\bibitem[3]{key-1} N.Nikolski, \textit{Condition Numbers of Large
Matrices and Analytic Capacities,} St. Petersburg Math. J., 17 (2006),
641-682.

\bibitem[4]{key-2} N.Nikolski, \textit{Operators, Function, and Systems:
an easy reading}, Vol.1. AMS, Providence, 2002.

\bibitem[5]{key-2} N.Nikolski, \textit{Treatise on the shift operator},
Springer-Verlag, Berlin etc., 1986.

\bibitem[6]{key-3} N. Nikolski, and S. A. Khrushchev, \textit{function
model and some problems in the spectral theory of functions}, Trudy
Mat. Inst. Steklov. 176, (1987), 97\textendash{}210 (Russian). English
transl.: Proc. Steklov Inst. Math. (1988), 101\textendash{}214.

\end{thebibliography}
\end{document}